\documentstyle[12pt, leqno]{article}

\begin{document}
\def\s{\subseteq}
\def\n{\noindent}
\def\se{\setminus}
\def\dia{\diamondsuit}
\def\la{\langle}
\def\ra{\rangle}

%--------------------------------------------------------------

\title{ Multiplicative Zagreb indices of k-trees}
\footnotetext{Partially support by College of Liberal Arts Summer Research Grant. }

\author{Shaohui Wang\footnote{  Corresponding author: S. Wang (e-mail: swang4@go.olemiss.edu), B. Wei (e-mail:  bwei@olemiss.edu). } , Bing Wei \\ 
\small\emph {Department of Mathematics, The University of Mississippi, University, MS 38677, USA}}
\date{}
\maketitle

Cite by(Shaohui Wang and Bing Wei, Multiplicative Zagreb indices of k-trees, Discrete Applied Mathematics, 180 (2015), 168-175, http://www.sciencedirect.com/science/
article/pii/S0166218X14003692.)

\begin{abstract}
Let $G$ be a graph with vetex set $V(G)$ and edge set $E(G)$.
 The first generalized multiplicative Zagreb index of $G$ is $\prod_{1,c}(G) = \prod_{v \in V(G)}d(v)^{c}$, for a real number $c > 0$, and the second multiplicative Zagreb index is $\prod_2(G) = \prod_{uv \in E(G)} d(u)d(v)$, where $d(u), d(v)$ are the degrees of the vertices of $u, v$. The multiplicative Zagreb indices have been the focus of considerable research in computational chemistry dating back to Narumi and Katayama in 1980s. In this paper, we generalize  Narumi-Katayama index and the first multicative index, where $c= 1, 2$, respectively,  and extend the results of Gutman to the generilized tree, the $k$-tree, where the results of Gutman are for $k = 1$. Additionally, we characterize the extremal graphs and determine the exact bounds of these indices of $k$-trees, which attain the lower and upper bounds.

\end{abstract}

\section{Introduction}
Throughout this paper $G = (V, E)$ is a  connected finite simple undirected graph with vertex set $V = V(G)$ and edge set $E =E(G)$.  Let $|G|$ or $|V|$ denote the cardinality of $V$.
For $S \subseteq V(G)$ and  $F \subseteq E(G)$,   we use $G[S]$ for the subgraph of $G$ induced by $S$, $G - S$ for the subgraph induced by $V(G) - S$ and $G - F$ for the subgraph of G obtained by deleting $F$.  Let $w(G - S)$ be the number of components of $G - S$, and  $S$ be a cut set if $w(G - S) \geq 2$.  For a vertex $v \in V(G)$, the neighborhood of  $v $ is the set $N(v) = N_G(v) = \{w \in V(G), vw \in E(G)\}$, and $d_G(v)$ (or briefly $d(v)$)  denotes the degree of $v$ with $ d_{G}(v) = |N(v)|$.  
 We use $G \cong H$ to denote that $G$ is isomorphic to $H$ and $G \not\cong H$ to denote  that $G$ is not isomorphic to $H$. Let $K_n, P_n$ and $S_n$ denote the clique, the path and the star on $n$ vertices, respectively. In particular, we say $K_n$ is a $k$-clique for $n = k$. 

In the 1980s, Narumi and Katayama [7] considered the product 
$$\begin{array}{rcl}
NK = \prod_{v \in V(G)} d(v)
\end{array}$$
which is  the "Narumi-Katayama index". And recently, Todeschini and Gutman et al [4, 10, 11]  studied the first and second multiplicative Zagreb indices defined as follow:
$$\begin{array}{rcl}
\prod_1(G) &=& \prod_{v \in V(G)}d(v)^{2},\\
\prod_2(G) &=& \prod_{uv \in E(G)} d(u)d(v).
\end{array}$$
Obviously, the first multiplicative Zagreb index is just the square of the NK index. Gutman [4] in 2011 characterized the multipilicative Zagreb indices for trees and determined the unique trees that obtained maximum and minimun values for $\prod_1(G)$ and $\prod_2(G)$, respectively.

\vskip 2mm {\bf Theorem 1 (Gutman 2011) } \emph{ Let $n \geq 5$ and  $T_n$ be any tree with $n$ vertices, then
$$\begin{array}{rcl}
(i) \prod_1(S_n) \leq \prod_1(T_n) \leq \prod_1(P_n);\\
 (ii)\prod_2(P_n) \leq \prod_2(T_n) \leq \prod_2(S_n).
\end{array}$$
  }

In this paper, we consider the first generalized multiplicative Zagreb index defined in (1) below and the second multiplicative Zagreb index: 
for any real number $c>0$,
$$\begin{array}{rcl}
(1)\prod_{1,c}(G) &=& \prod_{v \in V(G)}d(v)^{c};\\
(2)\prod_2(G) &=& \prod_{uv \in E(G)} d(u)d(v).
\end{array}$$
Eventually, for $c = 1, 2$,  (1) is just the NK index and the first multuplicative Zagreb, respectively. For (2), it is easy to see that $\prod_2(G) = \prod_{v\in V(G)}d(v)^{d(v)}$. Also we will find the bounds of the values of   $\prod_{1,c}(G)$, $\prod_{2}(G)$ for $k$-trees, respectively, and determine the extremal graphs which attain the bounds. Our main results are as follows:

{\bf Theorem 2 }  \emph 
Let $T_n^k$ be a $k$-tree on $n \geq k$ vertices, then
 $$\begin{array}{rcl}
  \prod_{1,c}(S_{k,n-k}) \leq \prod_{1,c} (T_n^k) \leq \prod_{1,c}(P_n^k),\end{array}$$
the left-side  and the right-side equalities are reached if and only if $T_n^k \cong S_{k,n-k}$ and $T_n^k \cong P_n^k$, respectively.

\vskip 2mm 
{\bf Theorem 3 }  \emph 
Let $T_n^k$ be a $k$-tree on $n \geq k$ vertices, then
 $$\begin{array}{rcl}
  \prod_2(P_n^k) \leq \prod_2 (T_n^k) \leq \prod_2(S_{k,n-k}),\end{array}$$
 the left-side  and the right-side equalities are reached if and only if $T_n^k \cong P_n^k$ and $T_n^k \cong S_{k,n-k}$, respectively.

\section{Preliminary}
It is commonly known that the class of $k$-trees is an important subclass of trangular graphs. Harry and Plamer [5] first introduced the $2$-tree in 1968, which is showed to be maximal outerplanar graphs in [3, 6]. Beineke and Pippert [1] gave the definition of $k$-trees in 1969. Relating to $k$-trees, there are many interesting applications to the study of a computational complexity and the intersection between graph theory and chemistry [2, 9]. We will just give some notations and definitions below.

\vskip 2mm {\bf Notation 1. } \emph{Let $[a, b]$ be the set of all the integers between $a$ and $b$ with $a \leq b$ including $a, b$, where $a, b$ are integers. Also, let $(a, b] = [a, b] -\{a\}$ and $[a, b) = [a, b] -\{b\}$. In particular, $[a, b] =\phi$ for $a >b$.
}
\vskip 2mm {\bf Notation 2. } \emph{
For ang integer p, if $p \ge 0$, we denote $x_{max\{0, p\}} = x_p$; If $p<0$, we say  $x_{max\{0, p\}} $ does not exist.
}
 
 \vskip 2mm {\bf Definition 1. } \emph{The $k$-tree, denoted by $T_n^k$, for  positive integers $n, k$ with $n \geq k$, is defined recursively as follows: The smallest $k$-tree is the $k$-clique $K_k$. If $G$ is a $k$-tree with $n \geq k$ vertices and a new vertex $v$ of degree $k$ is added and joined to the vertices of a $k$-clique in $G$, then the obtained graph is a $k$-tree with $n + 1$ vertices.
}

 \vskip 2mm {\bf Definition 2. } \emph{The $k$-path, denoted by $P_n^k$, for positive integers $n, k$ with $n \geq k$, is defined as follows:  Starting with a $k$-clique $G[\{v_1, v_2 ... v_k\}]$. For $i \in [k + 1, n]$, the vetex $v_i$ is adjacent to vertices $\{v_{i - 1}, v_{i - 2} ... v_{i - k} \}$ only.
}

\vskip 2mm {\bf Definition 3. } \emph{The $k$-star, denoted by $S_{k, n-k}$, for positive integers $n, k$ with $n \geq k$, is defined as follows:  Starting with a $k$-clique $G[\{v_1, v_2 ... v_k\}]$ and an independent set $S$ with $|S| = n - k$. For $i \in [k + 1, n]$, the vetex $v_i$ is adjacent to vertices $\{v_{1}, v_{ 2} ... v_{k} \}$ only.
}

\vskip 2mm {\bf Definition 4. } \emph{A vertex $v \in V(T_n^k)$ is called a $k$-simplicial vertex if $v$ is a vertex of degree $k$ whose neighbors form a $k$-clique of $T_n^k$. Let $S_1(T_n^k)$ be the set of all $k$-simplicial vertices of $T_n^k$, for $n \geq k+2$, and set $S_1(K_k) = \phi, S_1(K_{k+1}) = \{v\} $, where $v$ is any vetex of $K_{k+1}$. Let $G =G_0, G_i = G_{i - 1} - v_i$, where $v_i$ is a $k$-simplicial vetex of $G_{i - 1}$, then $\{v_1, v_2 ... v_n \}$ is called a simplicial elimination ordering of the $n$-vertex graph $G$.
}

\vskip 2mm {\bf Definition 5. } \emph{If   $w(G - S) \leq 2$ for   any $k$-clique $G[S]$ of $T_n^k$, we say $T_n^k$ is a hyper pendent edge; If there exists a $k$-clique $G[S]$ with $w(G - S) \geq 3$,   let $C$ be a component of $T_n^k - S$ and contain a unique vertex belonging to $S_1(G)$, then we say that $G[V(S) \cup V(C)]$ is a hyper pendent edge of $T_n^k$, denoted by $\mathcal{P}$. In particular, a $k$-path is a hyper pendent edge.
}

Moreover, let $G[\{v_1, v_2 ... v_k\}]$ denote the initial $k$-clique, then just by the definition of $k$-trees, one can get 
\vskip 2mm
{\bf Fact 1. }  \emph  {For the $k$-star, the degree of vertex $v_i$ can be characterized as follows: $d(v_i) = n - k$, for $i \in [1, k]$; $d(v_i) = k$, for $i \in [k+1, n]$. }
\vskip 2mm
{\bf  Fact 2. }  \emph  {For the $k$-path, the degree of vertex $v_i$ can be characterized as follows:\\
 (1) If $4 \leq n \leq 2k$, $d(v_i) = k + i -1 $, for $i \in [1, n - k -1]$; $d(v_i) = n - 1$, for $i \in [n - k, k + 1]$; $d(v_i) = k + n - i$, for $i \in [k + 2, n]$.\\
(2) If $n \geq 2k + 1$, $d(v_i) = k + i - 1$, for $i \in [1, k]$; $d(v_i) = 2k$, for $i \in [k+1, n-k]$; $d(v_i) = k + n -i$, for $i \in [n - k + 1, n]$.}

Easily verified through induction by using the above obseaverations, one can deduce the first generalized multiplicative Zagreb indices and second multiplicative Zagreb indices of the $k$-path and $k$-star as follows.

\vskip 2mm {\bf  Fact 3. }  \emph {
Let $S_{k,n-k}$ be a $k$-star on $n \geq k + 1$ vertices, then \\
$(1)  \prod_{1,c} (S_{k,n-k}) = (n - k)^{ck} k^{c(n-k)};$\\
$(2)  \prod_2 (S_{k,n-k}) = (n - k)^{k(n - k)} k^{k(n - k)}.$}

\vskip 2mm {\bf Fact 4. }  \emph {
Let $P_n^k$ be a $k$-path on $n \geq k + 1$ vertices, then \\
$
(1 . 1)  \prod_{1,c} (P_n^k) = (n - 1)^c \prod_{i = k}^{n - 2} i^{2c} $, if $n \in [k +1, 2k]$;\\
$
(1 . 2)  \prod_{1,c} (P_n^k) = (2k)^{c(n - 2k)} \prod_{i = k}^{2k - 1} i^{2c}$, if $n \geq 2k + 1;$\\
$
(2 . 1)  \prod_2 (P_n^k) = (n - 1)^{n - 1} \prod_{i = k}^{n - 2} i^{2i} $, if $n \in [k +1, 2k]$;\\
$
(2 . 2)  \prod_2 (P_n^k) = (2k)^{2k(n-2k)} \prod_{i = k}^{2k - 1} i^{2i}$, if $n \geq 2k + 1$.}

By considering the derivatives of the following functions, one can get

\vskip 2mm {\bf  Fact 5. }  \emph {
The function $f(x) =\displaystyle{\frac {x}{x+m}}$ is strictly increasing for  $x \in [0, \infty)$, where $m$ is a positive integer.
}

\vskip 2mm {\bf  Fact 6. }  \emph {
The function $f(x) =\displaystyle{\frac {x^x}{(x+m)^{x+m}}}$  is strictly decreasing for  $x \in [0, \infty)$, where $m$ is a positive integer.
}

\section{Main proofs }
Firstly, we give some lemmas that are critical in the proof of our main results.
 \vskip 2mm {\bf Lemma 1 } \emph{    For any $k$-tree $G \not\cong S_{k,n-k}$,  let $u \in S_2$, $N(u) \cap S_1 =\{v_1, v_2... v_s\}$, where $s \geq 1$ is an integer, then\\
(1) For any $i$ with $1 \leq i \leq s$, there exists a vertex $v \in N(u) - \{v_1, v_2 ... v_s\}$ of degree at least $k$ in $G[V(G) - \{v_1, v_2 ... v_s\}]$ such that $vv_i \notin E(G)$.      \\
(2)There exists a $k$-tree $G^*$ such that $\prod_{1,c}(G^*) < \prod_{1,c}(G)$ and  $\prod_2(G^*) > \prod_2(G)$.
  }

{\bf Proof.  } \emph
For (1),  let $G' = G[V(G) - \{v_1, v_2 ... v_s\}]$ and $S = N(u) - \{v_1, v_2 ... v_s\}$,   we obtain that $d_{G'}(u) = |S| = k$ and
$G[S]$ is a $k$-clique by $u \in S_2$. 
Since $G \not\cong S_n^k$, $d_{G'}(v) \ge k$ for all $v \in S$.  And by the facts that $N(v_i) \subseteq (N(u) - \{v_1, v_2 ... v_s\}) \cup \{u\}$ with $|N(v_i)| = k$ and $| (N(u) - \{v_1, v_2 ... v_s\}) \cup \{u\}| = k+1$, we have   for any $i \in [1,s]$, there exists a vertex $v \in S$ such that $vv_i \notin E(G)$.

For (2), choose $v_1$ and by (1) there exists a vertex $v \in N(u) - \{v_1, v_2 ... v_s\}$ with $d_{G'}(v) \geq k $ such that $vv_1 \notin E(G)$. If $d_{G'}(v) = k $, and by $uv \in E(G')$, we obtain $G'$ is a $k+1$-clique. Let  $x \in S$ be the vertex such that $d(x) = min_{v\in S} \{d(v)\}$, and let $v_t$ be the vertex such that $v_tx \in E(G)$, $v_ty \notin E(G)$ for some $t \in [1,s]$ and $y \in S$, that is, $ d(x) -1 < d(y)$. Construct a new graph $G^*$ such that $V(G^*) = V(G)$, and $E(G^*) = E(G) - \{v_tx\} + \{v_ty\}$.  Denote $G_0 = G[V(G) - \{x,y\}]$,  since $  d(x) -1 < d(y)$, and by the definition of $\prod_{1,c}(G)$,  $\prod_{2}(G)$ and  Fact 5, Fact 6,  we have

 $$\begin{array}{rcl}
 \displaystyle{ \frac {\prod_{1,c}(G)}{\prod_{1,c}(G^*)}} &=& 
 \displaystyle{\frac { [\prod_{w \in V(G_0)} d(w)^c] d(y)^c d(x)^c}{[\prod_{w \in V(G_0)} d(w)^c ] [d(y)+1]^c[d(x)-1]^c}}\\ 
&=&  \displaystyle{\frac {\displaystyle{d(y)^cd(x)^c}}{\displaystyle{[d(y) + 1]^c[d(x) -1]^c}} } \\
&=&  \displaystyle{\frac { \displaystyle{\frac {d(y)^c}{[d(y) + 1]^c}}} {\displaystyle{\frac{[d(x) - 1]^c} {d(x)^c}} }}\\
&>& 1.
\end{array}$$

Also, 

$$\begin{array}{rcl}  \displaystyle{
\frac {\prod_2(G)}{\prod_2(G^*)}} 
&=&   \displaystyle{ \frac { [\prod_{w \in V(G_0)}d(w)^{d(w)}] d(y)^{d(y)}d(x)^{d(x)}}{[\prod_{w \in V(G_0)}d(w)^{d(w)} ] [d(y)+1]^{d(y)+1}[d(x)-1]^{d(x)-1}}}  \\
&=&  \displaystyle{ \frac {d(y)^{d(y)}d(x)^{d(x)}}{[d(y)+1]^{d(y)+1}[d(x)-1]^{d(x)-1}}}\\
&=&   \displaystyle{\frac{\displaystyle{[\frac{d(y)^{d(y)}}{[d(y)+1]^{d(y)+1}}]}} { \displaystyle{[\frac {[d(x)-1]^{d(x)-1}}{d(x)^{d(x)}}]}} }\\
&<& 1.\\
\end{array}$$
Thus, we find that the $k$-tree $G^*$ satisfies $\prod_{1,c}(G^*) < \prod_{1,c}(G)$ and  $\prod_2(G^*) > \prod_2(G)$, we are done.

If $d_{G'}(v) \geq k+1 $,  reorder the subindices of $\{v_1, v_2 ... v_s\}$ such that  $vv_i \notin E(G)$ with $i \in [1, s_1]$, where $s_1 \leq s$, and
 by the fact that $G[N(u) - \{v_1, v_2 ... v_s\}]$ is a $k$-clique,  we have $d(u) = k+s$ and $d(v) \geq k+1+s-s_1 $, that is, $ d(v) \geq d(u)  -s_1 +1$.
Construct a new graph $G^*$ such that $V(G^*) = V(G)$, and $E(G^*) = E(G) - \{uv_i\} + \{vv_i\}$, for all $i \in [1,s_1]$.
Since $G[N(u)-\{v_1,v_2 ... v_s\}+ \{u\}]$ is a $k+1$-clique, and for any $i$, $N(v_i) \subseteq N_{G-\{v_1,v_2 ... v_s\}}(u) \cup \{u\}$, then  $G^*$ is a $k$-tree. Denote $G_0 = G[V(G) - \{u, v\}]$,  since $ d(v) \geq d(u)  -s_1 +1$, and by the definition of $\prod_{1,c}(G)$,  $\prod_{2}(G)$ and  Fact 5, Fact 6,  we have

 $$\begin{array}{rcl}
 \displaystyle{ \frac {\prod_{1,c}(G)}{\prod_{1,c}(G^*)}} &=& 
 \displaystyle{\frac { [\prod_{w \in V(G_0)} d(w)^c] d(v)^c d(u)^c}{[\prod_{w \in V(G_0)} d(w)^c ] [d(v)+s_1]^c[d(u)-s_1]^c}}\\ 
&=&  \displaystyle{\frac {\displaystyle{d(v)^cd(u)^c}}{\displaystyle{[d(v) + s_1]^c[d(u) - s_1]^c}} } \\
&=&  \displaystyle{\frac { \displaystyle{[\frac {d(v)^c}{[d(v) + s_1]^c}]}} {\displaystyle{[\frac{[d(u) - s_1]^c} {d(u)^c}]} }}\\
&>& 1.
\end{array}$$

Also, 

$$\begin{array}{rcl}  \displaystyle{
\frac {\prod_2(G)}{\prod_2(G^*)}} 
&=&   \displaystyle{ \frac { [\prod_{w \in V(G_0)}d(w)^{d(w)}] d(v)^{d(v)}d(u)^{d(u)}}{[\prod_{w \in V(G_0)}d(w)^{d(w)}]  [d(v)+s_1]^{d(v)+s_1}[d(u)-s_1]^{d(u)-s_1}}}  \\
&=&  \displaystyle{ \frac {d(v)^{d(v)}d(u)^{d(u)}}{[d(v)+s_1]^{d(v)+s_1}[d(u)-s_1]^{d(u)-s_1}}}\\
&=&   \displaystyle{\frac{\displaystyle{[\frac{d(v)^{d(v)}}{[d(v)+s_1]^{d(v)+s_1}}]}} { \displaystyle{[\frac {[d(u)-s_1]^{d(u)-s_1}}{d(u)^{d(u)}}]}} }\\
&<& 1.\\
\end{array}$$

Hence, we find that the $k$-tree $G^*$ satisfies $\prod_{1,c}(G^*) < \prod_{1,c}(G)$ and  $\prod_2(G^*) > \prod_2(G)$, we are done.
$\hfill\Box$

\vskip 2mm {\bf Lemma 2 } \emph{Let $G$ be a $k$-tree, if either  $\prod_{1,c}(G)$ attains the maximal or  $\prod_2(G)$ attains the minimal,  then every hyper pendent edge is a k-path. 
}

{\bf Proof.  } \emph
Let $\mathcal{P}$$ = G[V(S) \cup V(C)]$ be a hyper pendent edge, where $G[S]=G[\{x_1, x_2 ... x_k\}]$ is a cut $k$-clique and  $V(C) = \{u_1, u_2 ... u_p\}$ with $p$ is a positive ingeter such that  $u_1$ is the only vertex of  $\mathcal{P}$ in $S_1(G)$ and  for $i \in [1, p-1]$, $u_{i}$ is the vertex added following  by $u_{i+1}$ through the process of Definition 1. 

\vskip 2mm {\bf  Fact 7. }  \emph {
For any hyper pendent edge  $\mathcal{P}$$ = G[V(S) \cup V(C)]$ as represented above,   $\{u_1, u_2 ... u_{p}\}$ is a simplicial elimination ordering of $\mathcal{P}$.}

{Proof.} By contradiction, assume that  $\{u_1, u_2 ... u_{p}\}$ is not a simplicial elimination ordering of $\mathcal{P}$. Let $u_t$ be the first vertex from $u_1$ to $u_p$ such that $\{u_t, u_{t+1}\} \in S_t$ for $t \in [2, p-1]$, then  $u_{t}u_{t+1} \notin E(G)$ and $\{u_{t}, u_{t+1}\}$ can not be in some $k$-cliques. And by Definition 1, there must be at least two vertices that belongs to $S_1$ in $V(C)$, a contradiction. $\Box$  

By Fact 7, we know $\{u_1, u_2 ... u_{p}\}$ is a simplicial elimination ordering of $\mathcal{P}$.  For $p \leq 2$,  $\mathcal{P}$ is a $k$-path by Definition 2; For $p \geq 3$, if  $\mathcal{P}$ is a $k$-path, then we are done. Otherwise, let $u_s$ be the first vertex from $u_p$ to $u_1$ such that  $G[V(S) \cup \{u_{p}, u_{p-1} ... u_{s+1}, u_{s}\}]$ is not a $k$-path. Since $G[V(S) \cup \{u_{p}, u_{p-1} ... u_{s+1}\}]$ is a $k$-path,  for each $i\in[s+1, p]$, let $N_{G-\{u_1, u_2 ... u_{i-1}\}}(u_i) = \{u_{i+1}, u_{i+2} ... ...  u_{min\{p, i+k\}}, x_1, x_2 ... x_{ max\{0, k-p+i\}}\}$, and by Definition 2 and the symmetry of $G[S]$, we have $|N(u_s) \cap \{u_{s+1}, u_{s+2} ... u_{min\{p, s+k\}}\}| = min\{p-s-1,k-1\}$, where $1 \leq s \leq p - 1$.

 For $p \leq k+s$, 
suppose that $u_t$ is the vertex such that $u_t \notin N(u_s)$ with $s+2 \leq t \leq p$,  let $N_{G-\{u_1, u_2 ... u_{s-1}\}}(u_s) = \{u_{s+1}, u_{s+2} ... u_{t-1}, u_{t+1} ...  u_{p}, x_1, x_2 ... x_{ k - p + s+1}\} $, and let $|N(x_{k-p+s+1}) \cap \{u_1, u_2 ... u_{s-1}\}| = m $ for $m \in [0, s-1]$. By Definition 2, we have $u_tu_i \notin E(G)$ for $i \in [1, s]$, and then $d(u_{t}) = k + t-s-1$ and $d(x_{k-p+s+1}) > k + p -s +m -1 $. Now construct a new graph $G^*$ such that $V(G^*) = V(G),  E(G^*) = E(G) - \{u_sx_{k-p+s+1}, u_ix_{k-p+s+1}\} + \{u_su_{t}, u_iu_{t}\}$ with $i \in [0, m]$, then $G^*$ is a $k$-tree. Since $t \leq p$, we have $d(x_{k-p+s+1}) > d(u_{t}) +m +1$,  and by the definition of $\prod_{1,c}(G)$, $\prod_2(G)$ and Fact 5, Fact 6,  we get

$$\begin{array}{rcl}
\frac {\prod_{1,c}(G)}{\prod_{1,c}(G^*)} &=& \displaystyle{ \frac {d(u_t)^cd(x_{k-p+s+1})^c}{[d(u_t) +m+1]^c[d(x_{k-p+s+1})-m-1] ^c} } \\
  &=&  \displaystyle{\frac {[\displaystyle{\frac{d(u_t)} {d(u_t)+m+1}}]^c}{[\displaystyle{\frac{d(x_{k-p+s+1}) - m-1}{d(x_{k-p+s+1})}}]^c}} \\
&< & 1,
\end{array}$$
$$\begin{array}{rcl}
\displaystyle{ \frac {\prod_{2}(G)}{\prod_{2}(G^*)}} &=& \displaystyle{ \frac {d(u_t)^{d(u_t)}d(x_{k-p+s+1})^{d(x_{k-p+s+1})}}{[d(u_t) + m+1]^{d(u_t) +m+1}[d(x_{k-p+s+1}) - m-1]^{d(x_{k-p+s+1}) - m-1}} } \\
  &=&  \displaystyle{ \frac {\displaystyle{\frac {d(u_t)^{d(u_t)}}{[d(u_t)+m+1]^{d(u_t)+m+1}}}} {\displaystyle{\frac{[d((x_{k-p+s+1}) - m-1]^{d((x_{k-p+s+1}) - m-1}}{d(x_{k-p+s+1})^{d(x_{k-p+s+1})}}}}} \\
&> & 1.
\end{array}$$

Thus,  $\prod_{1,c}(G^*) > \prod_{1,c}(G)$ and $\prod_{2}(G^*) < \prod_{2}(G)$, a contradiction. 

 For $p \geq k+s +1$, let $|N(u_{k+s+1}) \cap \{u_1, u_2 ... u_{s-1}\}| = m$ for $m \in [0, s-1]$.  Since $G[V(S) \cup \{u_{p}, u_{p-1} ... u_{s+1}\}]$ is a $k$-path, we have $G[\{u_{s+1}, u_{s+2} ... u_{s+k+1}\}]$ is a $k+1$-clique. Suppose that $u_t$ is the vertex such that $u_t \notin N(u_s)$ with $s+2 \leq t \leq s+k$, let
 $N_{G- \{u_1, u_2 ... u_{s-1}\}}(u_s) = \{u_{s+1}, u_{s+2} ... u_{t-1}, u_{t+1} ... u_{s+k+1}\}$. Now we construct a new graph  $G^*$ such that $V(G^*) = V(G)$, $E(G^*) = E(G) - \{u_su_{k+s+1}, u_iu_{k+s+1}\} + \{u_su_{t}, u_iu_t\}$ for $i \in [0, m]$,  then $G^*$ is a $k$-tree and $d(u_{k+s+1}) = 2k+ m$, $d(u_{t})= k+t-s-1$.  Since $t \leq s+k$, we have $d(u_{k+s+1}) > d(u_{t})+m+1$, and by the definition of $\prod_{1,c}(G)$, $\prod_2(G)$ and Fact 5, Fact 6,  we get

$$\begin{array}{rcl}
\frac {\prod_{1,c}(G)}{\prod_{1,c}(G^*)} &=& \displaystyle{ \frac {d(u_t)^cd(u_{k+s+1})^c}{[d(u_t) +m+1]^c[d(u_{k+s+1})-m-1] ^c} } \\
  &=&  \displaystyle{\frac{[\displaystyle{\frac {d(u_t)}{d(u_t)+m+1}]^c}} {[\displaystyle{\frac{d(u_{k+s+1})-m - 1}{d(u_{k+s+1})}]^c}}} \\
&< & 1,
\end{array}$$
$$\begin{array}{rcl}
\displaystyle{ \frac {\prod_{2}(G)}{\prod_{2}(G^*)}} &=& \displaystyle{ \frac {d(u_t)^{d(u_t)}d(u_{k+s+1})^{d(u_{k+s+1})}}{[d(u_t) +m+ 1]^{d(u_t)+m +1}[d(u_{k+s+1})-m - 1]^{d(u_{k+s+1})-m - 1}} } \\
  &=&  \displaystyle{ \frac {\displaystyle{\frac{d(u_t)^{d(u_t)}} {[d(u_t)+m+1]^{d(u_t)+m+1}}}}{\displaystyle{\frac{[d((u_{k+s+1})-m - 1]^{d(u_{k+s+1}) -m- 1}}{d(u_{k+s+1})^{d(u_{k+s+1})}}}}} \\
&> & 1.
\end{array}$$

Thus,  $\prod_{1,c}(G^*) > \prod_{1,c}(G)$ and $\prod_{2}(G^*) < \prod_{2}(G)$, a contradiction. 
Hence, for any $s \in [1, p]$ $N_{G-\{u_1, u_2 ... u_{s-1}\}}(u_s) = \{u_{s+1}, u_{s+2} ... u_{min \{p, k+s\}}, x_1, x_2 ... x_{max \{0, k - p + s\}}\} $, that is,  $\mathcal{P}$ is a $k$-tree. $\hfill\Box$

\vskip 2mm {\bf Lemma 3 } \emph{    Let $G$ be a  $k$-tree, if either $\prod_{1,c}(G)$ attains the maximal or  $\prod_2(G)$ attains the minimal, then $|S_1(G)| = 2$.  
 }

{\bf Proof.  } \emph 
We know that $|S_1(G)| \geq 2$ for $n \geq k + 2$, and by Lemma 2, every hyper pendent edge is a $k$-path  for $\prod_{1,c}(G)$ to attain the maximal or  $\prod_2(G)$ to attain the minimal. If $|S_1(G)| = 2$, we are done; Suppose that $|S_1(G)| \geq 3$, it suffices to prove that there exists a graph $G'$ such that $|S_1(G')| = |S_1(G)| - 1$ with $ \prod_{1,c}(G') > \prod_{1,c}(G)$ and  $\prod_2(G') < \prod_2(G)$. 

\vskip 2mm {\bf  Fact 8. }  \emph {
 For any $k$-tree $G$ satisfying the conditions of Lemma 3, if $|S_1(G)| \geq 3$, then there exists a $k$-clique $G[S]$ such that $w(G-S) \geq 3$.}

{Proof.}  We will proceed by induction on $n = |G|$. For $n= k+3$, it's trivial; For $n \geq k+4$, assume that the fact is true for the $k$-tree $G$ with  $n < k+p$, and consider $n = k+p$. 
 If $|S_1(G)| \ge 4$, choose any vertex $v \in S_1(G)$, or $|S_1(G)| = 3$ and $|S_2(G)| \ge 2$, choose the vertex $v \in S_1(G)$ such that $ N(w) \cap S_1(G) = \{v\}$ for some $w \in S_2(G)$, then construct a new graph $G'$ such that $G' = G - v$. Since $S_2(G)$ is an dependent set and $G[N(v)]$ is a $k$-clique for any $v \in S_1(G)$, we obtain $|S_1(G')| \ge 3$. By the induction hypothesis, there exists a $k$-clique $G[S]$ in $G'$ such that $w(G' - S) \ge 3$.
Thus,  by adding back $v$, $G[S]$ is still a $k$-clique in $G$ and $w(G-S) \geq 3$, we are done. Next, we only consider  $|S_1(G)| = 3$ and $|S_2(G)| = 1$.

 Let $S_1(G) = \{v_1, v_2, v_3\}$ and $G_0 = G-\{v_1, v_2, v_3\}$,  by Definition 4, we have $G_0$  is a $k+1$-clique, denoted $G[\{x_1, x_2 ... x_{k+1}\}]$. If there exists $N(v_i) = N(v_j)$,  for some $i, j \in [1, 3]$ with $i\neq j$, and take $S = N(v_i)$, then $w(G- S) \ge 3$, we are done;  If $N(v_i) \neq N(v_j)$,  for any $i, j \in [1, 3]$ with $i \neq j$, then reorder the index of $x_i$ such that $N(v_1) = \{x_1, x_2 ... x_k\}$, $N(v_2) = \{x_2, x_3 ... x_{k+1}\}$ and $N(v_3) = \{x_1, x_3 ... x_{k+1}\}$. Construct a new graph $G^*$ such that $V(G^*) = V(G)$, $E(G^*) = E(G) - \{v_1x_1\} +\{v_1v_2\}$, then $G^*$ is still a $k$-tree and $d_G(x_1) = k+2$,  $d_{G^*}(x_1) = k+1$, $d_G(v_1) =d_G(v_2) =k$ and  $d_{G^*}(v_2) = k+1$.  By the definition of $\prod_{1,c}(G)$, $\prod_2(G)$ and Fact 6, we have 
$$\begin{array}{rcl}
\displaystyle{\frac {\prod_{1,c}(G)}{\prod_{1,c}(G^*)}} &=& \displaystyle{ \frac {d(v_2)^cd(x_{1})^c}{[d(v_2) +1 ]^c[d(x_{1})-1] ^c} } \\
  &=&  \displaystyle{[\frac {k(k+2)}{(k+1)^2}]^c} \\
&< & 1,
\end{array}$$
$$\begin{array}{rcl}
\displaystyle{ \frac {\prod_{2}(G)}{\prod_{2}(G^*)}} &=& \displaystyle{ \frac {d(v_2)^{d(v_2)}d(x_{1})^{d(x_{1})}}{[d(v_2) +1]^{d(v_2)+1 }[d(x_{1})-1]^{d(x_{1})-1 }} } \\
  &=&  \displaystyle{ \frac {(k+2)^{k+2}k^k}{(k+1)^{2(k+1)}}} \\
&=& \displaystyle{ \frac {\displaystyle{[\frac {k^k}{(k+1)^{k+1}}]} }{\displaystyle{[\frac {(k+1)^{k+1}}{(k+2)^{k+2}}]} }  } \\
&> & 1.
\end{array}$$

Thus, we find a graph $G^*$ with $\prod_{1,c}(G^*) > \prod_{1,c}(G)$ and $\prod_{2}(G^*) < \prod_{2}(G)$, a contradiction with that $\prod_{1,c}(G)$ attains the maximal or  $\prod_2(G)$ attains the minimal, we are done.
   $\Box$

Choose a $k$-clique $G[S]$ with $w(G - S) \geq 3$ such that there are two components of $G - S$: $C_1, C_2$ with $|C_1|=p, |C_2| = q$ and $p + q$ being minimal, for $p \geq q$. Let $u_1\in V(C_1)$,  $v_1\in V(C_2)$ with $\{u_1, v_1\} \subseteq S_1(G)$. Let $N_{G- \{v_1,v_2 ... v_{i-1}\}}(v_i) = \{v_{i+1}, v_{i+2} ... v_{min\{k+1, q\}}, x_1, x_2 ... x_{max\{0, {k-q+i}\}}$, 
$N_{G-\{u_1,u_2 ... u_{j-1}\}}(u_j) = \{u_{j+1}, u_{j+2} ... u_{min\{k+1, p\}}, y_1,y_2 ... y_{max\{0, {k-p+i}\}}\}$ for $i \geq 1, j \geq 1$, where $\{v_1, v_2 ... v_q\}$ and $\{u_1, u_2 ... u_p\}$ are simplicial elimination orderings of $G[S\cup V(C_1)]$ and $G[S\cup V(C_2)]$, respectively.  We will prove Lemma 3 by induction on $q$.

$(1)$ If $q =1$, then  $d(v_1) = k$. Choose $x_t \in N(v_1)$, let $|N(x_t) \cap \{u_1, u_2 ... u_p\}| = m$ for $m \in [1, k]$,  we get $d(x_t) > k+1+m$ by $w(G -S) \geq 3$, and then $d(x_t) > d(v_1) + m +1$. Now construct a new graph  $G^*$ such that $V(G^*) = V(G), E(G^*) = E(G) - \{u_ix_t\} + \{u_iv_{1}\}$ for $i \in [1, m]$,  then $G^*$ is a $k$-tree and $|C_1| + |C_2| = p$ with $G[\{x_1, x_2 ... x_{t-1}, x_{t+1} ... x_k, v_1\}]$ is a $k$-clique in $G^*$.  Since   $d(x_{t}) > d(v_{1})+m+1$,  by the definition of $\prod_{1,c}(G)$, $\prod_2(G)$ and Fact 5, Fact 6, we have

$$\begin{array}{rcl}
\frac {\prod_{1,c}(G)}{\prod_{1,c}(G^*)} &=& \displaystyle{ \frac {d(v_1)^cd(x_{t})^c}{[d(v_1) +m ]^c[d(x_{t})-m] ^c} } \\
  &=&  \displaystyle{\frac {[\displaystyle{\frac {d(v_1)}{d(v_1)+m}]^c}}{[\displaystyle{\frac{d(x_{t})-m }{d(x_{t})}]^c}}} \\
&< & 1,
\end{array}$$
$$\begin{array}{rcl}
\displaystyle{ \frac {\prod_{2}(G)}{\prod_{2}(G^*)}} &=& \displaystyle{ \frac {d(v_1)^{d(v_1)}d(x_{t})^{d(x_{t})}}{[d(v_1) +m]^{d(v_1)+m }[d(x_{t})-m]^{d(x_{t})-m }} } \\
  &=&  \displaystyle{ \frac {\displaystyle{\frac{d(v_1)^{d(v_1)}} {[d(v_1)+m ]^{d(v_1)+m }}}}{\displaystyle{\frac{[d((x_{t})-m ]^{d((x_{t}) -m}}{d(x_{t})^{d(x_{t})}}}}} \\
&> & 1.
\end{array}$$

Then,  $\prod_{1,c}(G^*) > \prod_{1,c}(G)$ and $\prod_{2}(G^*) < \prod_{2}(G)$. Thus,  let $G' = G^*$,  $|S_1(G')| = |S_1(G)| - 1$, $ \prod_{1,c}(G') > \prod_{1,c}(G)$ and  $\prod_2(G') < \prod_2(G)$, and we are done.

$(2)$ Assume that $q =s $, there exists a $k$-tree $G'$ such that $|S_1(G')| = |S_1(G)| -1$, $ \prod_{1,c}(G') > \prod_{1,c}(G)$,  $\prod_2(G') < \prod_2(G)$ and we consider $q = s+1 $.

If  $q \leq k$,  we have  $d(v_q) = k+q-1$ by the fact that $G[S\cup V(C_2)]$ is a $k$-path.
 Choose $x_t \in N(v_1)$,  we know $x_t \in N(v_i)$ for all $i \in [1, p]$ by  $G[S\cup V(C_2)]$ is a $k$-path. Let $|N(x_t) \cap \{u_1, u_2 ... u_p\}| = m$ for $m \in [1, k]$,  we have  $d(x_t) > k+q+m$ by $w(G - S) \geq 3$, and then $d(x_t) > d(v_q) + m +1$. Now construct a new graph  $G^*$ such that $V(G^*) = V(G), E(G^*) = E(G) - \{u_ix_t\} + \{u_iv_{q}\}$ for $i \in [1, m]$,  then $G^*$ is a $k$-tree and $|C_1| + |C_2| = p+q-1$ with $G[\{x_1, x_2 ... x_{t-1}, x_{t+1} ... x_k, v_q\}]$ is a $k$-clique in $G^*$.  Since   $d(x_{t}) > d(v_{q})+m+1$, by the definition of $\prod_{1,c}(G)$, $\prod_2(G)$ and Fact 5, Fact 6, we have

$$\begin{array}{rcl}
\frac {\prod_{1,c}(G)}{\prod_{1,c}(G^*)} &=& \displaystyle{ \frac {d(v_q)^cd(x_{t})^c}{[d(v_q) +m ]^c[d(x_{t})-m] ^c} } \\
  &=&  \displaystyle{\frac{[\displaystyle{\frac {d(v_q)}{d(v_q)+m}}]^c} {[\displaystyle{\frac{d(x_{t})-m }{d(x_{t})}}]^c}} \\
&< & 1,
\end{array}$$
$$\begin{array}{rcl}
\displaystyle{ \frac {\prod_{2}(G)}{\prod_{2}(G^*)}} &=& \displaystyle{ \frac {d(v_q)^{d(v_q)}d(x_{t})^{d(x_{t})}}{[d(v_q) +m]^{d(v_q)+m }[d(x_{t})-m]^{d(x_{t})-m }} } \\
  &=&  \displaystyle{ \frac {\displaystyle{\frac {d(v_q)^{d(v_q)}}{[d(v_q)+m ]^{d(v_q)+m }}}}{\displaystyle{\frac{[d((x_{t})-m ]^{d((x_{t}) -m}}{d(x_{t})^{d(x_{t})}}}}} \\
&> & 1.
\end{array}$$

Then, $ \prod_{1,c}(G) < \prod_{1,c}(G^*)$, $\prod_2(G) > \prod_2(G^*)$ and $q =s$ in $G^*$, then by the induction hypothesis, there exists a $k$-tree $G'$ such that $|S_1(G')| = |S_1(G)| -1$, $ \prod_{1,c}(G') > \prod_{1,c}(G)$ and  $\prod_2(G') < \prod_2(G)$, we are done.

 If $q \geq k+1$,   we have  $N(u_1) = \{u_2, u_3 ... u_{k+1}\}$, $N(v_1) = \{v_2, v_3 ... v_{k+1}\}$ by the facts that $p\geq q$ and $G[S\cup V(C_1)]$, $G[S\cup V(C_2)]$ are  $k$-paths. We construct a new graph $G^*$ such that  $V(G^*) = V(G), E(G^*) = E(G) - \{v_1v_i\} + \{u_jv_{1}\}$ for $i \in [2, k+1]$, $j \in [1, k]$. And by Fact 2 and the definition of $\prod_{1,c}(G)$, $\prod_2(G)$,
we obtain 
$$\begin{array}{rcl}$$
 \displaystyle{\frac {\prod_{1,c}(G)}{\prod_{1,c}(G^*)}} &=&
\displaystyle{ \frac { \prod_{i=2}^{k+1}d(v_i)^c \prod_{j=1}^{k}d(u_j)^c}{ \prod_{i=2}^{k+1}[d(v_i)-1]^c 
\prod_{j=1}^{k, }[d(u_j)+1]^c}}     \\
&=& 1,\\
 \displaystyle{\frac {\prod_2(G)}{\prod_2(G^*)}} &=&
\displaystyle{  \frac { \prod_{i=2}^{k+1}d(v_i)^{d(v_i)} \prod_{j=1}^{k}d(u_j)^{d(u_j)}}{ \prod_{i=2}^{k+1}[d(v_i)-1]^{d(v_i)-1} 
\prod_{j=1}^{k}[d(u_j)+1]^{d(u_j)+1}}
}\\
&=& 1.
\end{array}$$

Then, $ \prod_{1,c}(G) = \prod_{1,c}(G^*)$, $\prod_2(G) = \prod_2(G^*)$ and $q =s$ in $G^*$, then by the induction hypothesis, there exists a $k$-tree $G'$ such that $|S_1(G')| = |S_1(G)| -1$, $ \prod_{1,c}(G') > \prod_{1,c}(G)$ and  $\prod_2(G') < \prod_2(G)$, we are done.
 $\hfill\Box$

Now, we turn to prove the main results of the paper.

{\bf Proof of Theorem 2.  } \emph
For any $k$-tree $T_n^k$, if $|S_1(T_n^k)| = n - k$, then $T_n^k \cong S_{k,n-k}$, we are done.
And if $|S_1(T_n^k)| \leq n - k -1$, we can recursively use Lemma 1 to make $\prod_{1,c}(T_n^k)$ decreasing until $|S_1(T_n^k)| = n - k$. Thus, we have  $T_n^k \cong S_{k,n-k}$ for  $\prod_{1,c}(T_n^k)$ to arrive the minimal value.

By Lemma 2, if $\prod_{1,c}(T_n^k)$ get the maximal, then every hyper pendent edge is a $k$-path, and by Lemma 3, $|S_1(T_n^k)| = 2$, implying that  $T_n^k \cong P_n^k$ for  $\prod_{1,c}(T_n^k)$ to arrive the maximal value.    $\hfill\Box$

 \vskip 2mm 
{\bf Proof of Theorem 3.  } \emph
For any $k$-tree $T_n^k$, if $|S_1(T_n^k)| = n - k$, then $T_n^k \cong S_{k,n-k}$ ,we are done.
And if $|S_1(T_n^k)| \leq n - k -1$, we can recursively use Lemma 1 to make $\prod_2(T_n^k)$ increasing until $|S_1(T_n^k)| = n - k$, then we have $T_n^k \cong S_{k,n-k}$ for  $\prod_2(T_n^k)$ to arrive the maximal value.

By Lemma 2, if $\prod_2(T_n^k)$ get the minimal,  every hyper pendent edge is a $k$-path, and by Lemma 3, $|S_1(T_n^k)| = 2$. Then this $k$-tree is a $k$-path, that is, $T_n^k \cong P_n^k$ for  $\prod_2(T_n^k)$ to arrive the minimal value.  
 $\hfill\Box$

\end{document}